\documentclass[12pt]{amsart}

\usepackage{amssymb}
\usepackage{mathrsfs}

\def\F{\mathbb{F}}

\def\Z{\mathbb{Z}}

\newcommand{\pth}[1]{{#1}^{[p]}}

\def\su{{\subseteq}}
\def\la{{\langle}}
\def\ra{{\rangle}}

\def\w{{\omega}}

\def\dim{{\rm dim~}}

\def\ad{\mbox{ad~}}

\def\Proof{\noindent{\sl Proof.}\ }
\def\qed{{\hfill $\Box$ \medbreak}}

\newtheorem{defi}{Definition}[section]
\newtheorem{thm}[defi]{Theorem}
\newtheorem{lem}[defi]{Lemma}

\newtheorem{eg}[defi]{Example}
\newtheorem{prop}[defi]{Proposition}

\newtheorem*{thmm}{Main Theorem}

\begin{document}

\title[Non-Matrix PI enveloping algebras]{Non-Matrix Polynomial identity  enveloping algebras}
\author{\textsc{Hamid Usefi}}
\thanks{ Research was supported by  NSERC}
\date{\today}

\address{Department of Mathematics and Statistics,
Memorial University of Newfoundland,
St. John's, NL,
Canada, 
A1C 5S7}
\email{usefi@mun.ca}

\begin{abstract}
Let $L$ be a restricted Lie superalgebra with  its restricted enveloping algebra $u(L)$ over a field $\F$ of characteristic $p>2$.
A polynomial identity is called non-matrix if it is not satisfied by the algebra of $2\times 2$ matrices over $\F$.
We characterize $L$  when $u(L)$ satisfies a non-matrix polynomial identity.
\end{abstract}


  \maketitle

\section{Introduction}
A variety of associative algebras over a field $\F$ is called non-matrix if it does not contain $M_2(\F)$, the algebra of $2\times 2$ matrices over $\F$. A polynomial identity (PI) is called  non-matrix if $M_2(\F)$ does not satisfy this identity.
Latyshev in his attempt to solve the Specht problem proved that any non-matrix variety
generated by a finitely generated algebra over a field of characteristic zero is finitely based \cite{L77, L80}.
The complete solution of the Specht problem in the case of characteristic zero is given by Kemer \cite{Kemer2}.

Although several counterexamples are found for the Specht problem in the positive characteristic \cite{AK},  the development in this area has lead to some interesting results. Kemer in \cite{K96}  investigated the relation between PI-algebras and nil algebras and asked whether the Jacobson radical of a relatively free algebra of countable rank
over an infinite field of positive  characteristic is a nil ideal of bounded index. Amitsur had already proved in \cite{Am}   that the Jacobson radical of a relatively-free algebra of countable rank is nil and Samoilov in \cite{Sam} proved that the Jacobson radical of a relatively free algebra of countable rank over an infinite field of positive characteristic is a nilideal of bounded index. The non-matrix varieties have been further studied in \cite{BRT, MPR, R97}.

Enveloping algebras satisfying polynomial identities were first considered by Latyshev \cite{L63} by proving that
the universal enveloping algebra of a Lie algebra $L$  over a field of characteristic zero satisfies a PI if and only if $L$ is abelian. Latyshev's result was extended to positive characteristic by Bahturin \cite{B74}.
Passman \cite{P90} and Petrogradsky \cite{P91} considered the analogous problem for restricted Lie algebras and their envelopes.

Let $L=L_0\oplus L_1$ be a restricted Lie superalgebra with the bracket $(\,,)$.
 We denote the restricted enveloping algebra of $L$ by $u(L)$. All algebras in this paper are over a field $\F$ of characteristic $p>2$ unless otherwise stated. In case $p=3$ we add the axiom $((y,y),y)=0$, for every $y\in L_1$.
 This identity is necessary to embed $L$ in $u(L)$.
 Restricted Lie superalgebras whose enveloping algebras satisfy a polynomial identity have been characterized by Petrogradsky \cite{P92}. The purpose of this paper is to characterize restricted Lie superalgebras whose restricted enveloping algebras satisfy a non-matrix PI. Riley and Wilson considered similar conditions for restricted enveloping algebras and group algebras in \cite{RW99}. Recall that a subset $X\su L_0$ is called $p$-nilpotent if there exists an integer $s$ such that $x^{p^s}=0$, for every $x\in X$.
 Our main result is as follows.

\begin{thmm}\label{uL-non-matrix}
Let $L=L_0\oplus L_1$ be a restricted Lie superalgebra over a perfect field and denote by $M$  the  subspace spanned by all $y\in L_1$ such that $(y,y)$ is $p$-nilpotent.
The following statements are equivalent:
\begin{enumerate}
\item $u(L)$ satisfies a non-matrix PI.
\item $u(L)$ satisfies a PI, $(L_0, L_0)$ is $p$-nilpotent,  $\dim L_1/M\leq 1$,   $(M, L_1)$ is $p$-nilpotent, and  $(L_1, L_0)\su M$.
\item The commutator ideal of $u(L)$ is nil of bounded index.
\end{enumerate}
\end{thmm}
  Theorem \ref{Petrogradsky} below recalls Petrogradsky's characterization of when $u(L)$ satisfies an arbitrary PI strictly in terms of the underlying Lie superalgebra structure of $L$; this allows one to replace (2) with a similar such characterization (that is too cumbersome to state here).  Furthermore, we show that (2) implies (3) over any field.
However, given that $u(L)$ satisfies a non-matrix PI,  the restriction on the field is necessary to be able to show that $\dim L_1/M\leq 1$. In Section \ref{examples},  we show that over a non-perfect field there exists  a restricted Lie superalgebra $L=L_0\oplus L_1$ such that  $\dim L_1=2$, commutator ideal of $u(L)$ is nil of  index $2p$ and yet $(y, y)$ is not $p$-nilpotent, for every $y\in L_1$.
This is in complete contrast with the enveloping algebras of ordinary Lie superalgebras satisfying a non-matrix PI, see Theorem 1.2 of \cite{BRU} where a similar characterization does not require any restriction on the field.

\section{Preliminaries}
Unless otherwise stated, all algebras are over  a field $\F$ of characteristic $p>2$.
Let  $A=A_0\oplus A_1$ be a vector space decomposition of a non-associative algebra over $\F$. We
say that this is a $\mathbb{Z}_2$-grading of $A$ if $A_iA_j\su A_{i+j}$, for every
$i,j\in\Z_2$ with the understanding that the addition $i+j$ is mod 2.
The components $A_0$ and $A_1$ are called even and odd parts of $A$, respectively.
Note that $A_0$ is a subalgebra of $A$. One can associate a Lie super-bracket to $A$ by defining
$(x,y)=xy-(-1)^{ij}yx$ for every $x\in A_i$ and $y\in A_j$. If $A$ is associative, then
for any $x\in A_i$, $y\in A_j$ and $z\in A$ the following identities hold:
\begin{enumerate}
\item[(1)] $(x,y)=-(-1)^{ij}(y,x)$,
\item[(2)] $(x,(y,z))=((x,y),z)+(-1)^{ij}(y,(x,z)).$
\end{enumerate}
The above identities are the defining relations of Lie superalgebras. 
Furthermore, $A$ can be viewed as a Lie algebra by the usual Lie bracket
$[u,v]=uv-vu$.

If $L$ is a Lie superalgebra, we denote the bracket of $L$ by $( , )$.
The adjoint representation of $L$ is given by $\ad x : L\to L$, $\ad x(y) = (y, x )$, for all $x,y\in L$.
The notion of restricted Lie superalgebras can be easily formulated as follows:
\begin{defi}\label{def:res}
A Lie superalgebra $L = L_0\oplus L_1$ is called  \em
{restricted}, if there is a $p$th power map
$L_0\to L_0$, denoted as $\pth{}$, satisfying
\begin{enumerate}
\item[(a)] $\pth{(\alpha x)} = \alpha^p \pth x$, for all  $x \in L_0$ and $\alpha \in \F$,
\item[(b)] $(y, \pth x) = (y,_p x)=(\ad x)^p(y)$, for all $x \in L_0$
and $y \in L$,
\item[(c)] $\pth{(x + y)} = \pth x + \pth y + \sum_{i= 1}^{p-1}
s_i(x,y)$, for all $x, y \in L_0$ where $is_i$ is the
coefficient of ${\lambda}^{i-1}$ in $(\ad(\lambda x +y))^{p-1}(x)$.
\end{enumerate}
\end{defi}
For example, every $\mathbb{Z}_2$-graded associative algebra inherits a restricted Lie
superalgebra structure.

Let $L$ be a restricted Lie superalgebra. We denote the (restricted) enveloping algebra of $L$ by $u(L)$.
The augmentation ideal $\w(L)$ is the ideal of $u(L)$ generated by $L$.
The analogue of the PBW Theorem is as follows. We refer to \cite{BMPZ} for basic background.

\begin{thm}
Let $L=L_0\oplus L_1$ be a restricted  Lie superalgebra and let ${\mathcal B}$ be a totally ordered
basis for $L$ consisting of $\Z_2$-homogeneous elements. Then $u(L)$ has a basis
consisting of PBW monomials, that is, monomials of the form $x_{1}^{a_1}\ldots x_{s}^{a_s}$ where
$x_1<\cdots < x_s$ in ${\mathcal B}$, $0\leq a_i<p$ whenever $x_i\in L_0$, and $0\leq a_i\leq 1$ whenever $x_i\in L_1$.
\end{thm}

Since $L$ embeds into $u(L)$ we identify $\pth x$ with $x^p$, for every $x\in L_0$.
Note that $u(L)$ can be viewed as a Lie algebra via the Lie bracket $[x, y]=xy-yx$ and
if $x\in L_0$ and $y\in L$ then the bracket $(x,y)$ in $L$ is the same as the bracket $[x, y]$ in $u(L)$.
Let $H$ be a subalgebra of $L$. We denote by $H'$ the commutator subalgebra of $H$, that is $H'=(H, H)$.
For a subset $X\su  L$, we denote by $\la X\ra_p$ or $X_p$ the restricted ideal of $L$ generated by $X$. Also,
 $\la X\ra_{\F}$ denotes the subspace spanned by $X$.
 An element $x\in L_0$ is called \emph{$p$-nilpotent} if there exists some non-negative integer $t$ such that $x^{p^t}=0$. Also, recall that $X$ is said to be $p$-nil if every element $x\in X$ is $p$-nilpotent and  $X$ is $p$-nilpotent if there exists a positive integer $k$ such that $x^{p^k}=0$, for every $x\in X$.
 By an ideal of $L$ we always mean a restricted ideal, that is $I$ is an ideal of $L$ if
 $(I, L)\su I$ and $I_0$ is closed under the $p$-map.

Let  $B$ and $C$ be subspaces of $L$.
We denote by $(B, C)$ the subspace spanned by all commutators $(b, c)$, where $b\in B$ and $c\in C$.
The lower Lie central series of $L$ is defined  by setting
 $\gamma_1(L)=L$ and  $\gamma_n(L)=(\gamma_{n-1}(L), L)$, for every $n\geq 2$.
  Recall that $L$ is called nilpotent if $\gamma_n(L)=0$, for some $n$.
 The derived subalgebra of $L$ is defined by setting
$\delta_0(L)=L$ and $\delta_{i+1}(L)=(\delta_i(L), \delta_i(L))$, for every $i\geq 0$. Also, $L$ is called solvable if
$\delta_m(L)=0$, for some $m$, and the least of such $m$ is called the derived length of $L$.
Moreover, long commutators are left tapped, that is $(x,y,x)=((x,y), z)$.

Note that Engel's Theorem holds for Lie superalgebras, see \cite{sch}, for example.
\begin{thm}[Engel's Theorem]\label{Engel's Theorem}
Let $L$ be a finite-dimensional Lie superalgebra such that
$\text{ad } x$ is nilpotent, for every homogeneous element $x\in L$.
Then $L$ is nilpotent.
\end{thm}

The proof of the following lemma follows from Engel's Theorem and the fact that $(\text{ad } x)^2=\frac{1}{2} \text{ad } (x, x)$, for every $x\in L_1$.

\begin{lem}\label{L-res-nilpotent}
Let $L=L_0\oplus L_1$ be a finite-dimensional restricted Lie superalgebra.
If $L_0$ is  $p$-nil then $L$ is nilpotent.
\end{lem}

\begin{lem}\label{w(L)-nilpotent}
Let $L$ be a restricted Lie superalgebra. Then $\w(L)$ is associative nilpotent if and only if
$L$ is finite-dimensional  and $L_0$ is $p$-nil.
\end{lem}
\Proof The if part  follows from the PBW Theorem.
We prove the converse by induction on $\dim L$. By Lemma \ref{L-res-nilpotent}, $L$ is nilpotent and so
there exists a non-zero element $z$ in the center $Z(L)$ of $L$. Since $Z(L)$ is homogeneous we may assume that
either $z\in L_0$ or $z\in L_1$. If $z\in L_1$ then $z^2=(z,z)=0$. If $z\in L_0$ then
we can replace $z$ with its $p$-powers so that $z^p=0$. So in either case $z^p=0$ in $u(L)$.
Now consider $H=L/\la z\ra_p$. Then by induction hypothesis $\w(H)$ is nilpotent.
This means that $\w^m(L)\su \la z\ra_p u(L)$, for some $m$. It then follows that
$\w^{mp}(L)\su \la z^p\ra_p u(L)=0$, as required.
\qed

We shall use the following two results.

\begin{thm}[\cite{P92}]\label{Petrogradsky}
Let $L=L_0\oplus L_1$ be a restricted Lie superalgebra.
Then $u(L)$ satisfies a  PI if and only if there exist homogeneous restricted  ideals $B\su A\su L$ such that
\begin{enumerate}
\item $L/A$ and $B$ are both finite-dimensional.
\item $A'\su B$, $B'=0$.
\item The restricted Lie subalgebra  $B_0$ is $p$-nilpotent.
 \end{enumerate}
\end{thm}

\begin{prop}[\cite{RW99}]\label{ass-PI}
Let $R$ be an associative algebra that satisfies a non-matrix PI over a field of positive characteristic $p$. Then  there exists an integer $t$ such that $R$ satisfies the identity $([u, v]w)^{p^t}=0$.
\end{prop}

\section{Proofs}
\begin{prop}\label{M}
Let $L$ be a restricted Lie superalgebra over a perfect field $\F$. Let $M$ be the set consisting of  all $y\in L_1$ such that $(y,y)$ is $p$-nilpotent. If $u(L)$ satisfies a non-matrix PI then the following conditions are satisfied:
\begin{enumerate}
\item $M$ is a subspace of $L_1$ and $(L_1, L_0)\su M$;
\item $\dim L_1/M\leq 1$;
\item $(L_0, L_0)$  and $(M, L_1)$ are both $p$-nilpotent.
\end{enumerate}

\end{prop}
\Proof 
Let  $R=[u(L), u(L)] u(L)$. Proposition \ref{ass-PI} implies that $R$ is nil.
Let $y\in M$ and $x\in L_1$.
 Note that $(x, y)=2yx$ modulo $R$. So, $(y, x)^{2m}=(y, y)^m(x, x)^m$ modulo $R$. Since $(y, y)$ is $p$-nilpotent and $R$ is nil, we deduce that $(x, y)$ is $p$-nilpotent. So, $(M, L_1)$ is $p$-nil.
 Let $y_1, y_2\in M$ and set $y= \alpha y_1+\beta y_2$, where $\alpha, \beta\in \F$. Note that  
 $(y, y)^{p^i}=(y, \alpha y_1)^{p^i}+(y, \beta y_2)^{p^i}$ modulo $R$, for every $i$.
Since $(M, L_1)$ and $R$ are both nil, we deduce that $(y, y)$ is  $p$-nilpotent.
This proves that $M$ is a subspace.
Since $[L_1, L_0]\su R$,  we know by  Proposition \ref{ass-PI}  that $[L_1, L_0]$ is nil. Thus, if $y\in [L_1, L_0]$ then
$(y, y)=2y^2$ is $p$-nilpotent. Hence, $[L_1, L_0]\su M$
and this finishes the proof of (1).

Now, we show that $\dim L_1/M\leq 1$. Let $I$ be the subset of $L_0$ consisting of  $p$-nilpotent elements.
By Proposition \ref{ass-PI}, $I$ is a restricted ideal of $L_0$ and $(L_0, L_0)\su I$. We have also proved that  $ (M, L_1)\su I$. Note that $I+M$ is a restricted ideal of $L$.
Without loss of generality, we can replace $L$ with $L/ (I+M)$. So,  $(y, y)$ is not $p$-nilpotent, for every  $y\in L_1$. 
Let $y, z\in L_1$.
Since $(L_0, L)=0$,   we have
$$
[y, z]=-(y,z)+2yz, \quad [y, z]^2=(y, z)^2-(y,y)(z, z).
$$
By Proposition \ref{ass-PI}, there exists $m$ such that $[y, z]^{2p^m}=0$. So,  
\begin{align}\label{imprtant-equation}
(y, z)^{2p^m}=(y, y)^{p^m}(z, z)^{p^m}.
\end{align} 
Thus, by the PBW theorem, $(y, y)^{p^m}$ and $(z, z)^{p^m}$ must be linearly dependent.
So, $(z, z)^{p^m}=\beta (y, y)^{p^m}$, for some $\beta\in \F$. Equation \eqref{imprtant-equation} then implies that
\begin{align}\label{imprtant-equation-rev}
((y, z)^{p^m})^2=\beta ((y, y)^{p^m})^2.
\end{align} 
 Using the PBW Theorem again, we deduce that $(y, y)^{p^m}=\alpha (y, z)^{p^m}$, for some $\alpha \in \F$.
 Equation \eqref{imprtant-equation-rev} implies that $\beta \alpha^2=1$.
So, we must have $(y,y)^{p^m}=\alpha (y, z)^{p^m}$ and $(z, z)^{p^m}=\alpha^{-1} (y, z)^{p^m}$. Let $\gamma\in \F$ so that $\gamma^{p^m}=\alpha$. We have, $(y,y-\gamma  z)^{p^m}=0$ and $(z,z-\gamma^{-1}  y)^{p^m}=0$.
So, $(y-\gamma  z, y-\gamma  z)^{p^m}=0$ which implies that $y-\gamma  z=0$. Thus, $y$ and $z$ are dependent, as required.  This proves (2). In order to prove (3), it suffices to show that there exists an 
integer $m$ such that  $((L_0, L_0)+(M, L_1))^{p^m}=0$. Note that, by Proposition \ref{ass-PI}, there exists an integer $t$ such that $[u, v]^{p^t}=0$, for all $u, v\in u(L)$.
By Theorem \ref{Petrogradsky}, there exists a homogeneous ideal $A$ of $L$ of finite codimension
such that  $B=\la A'\ra_p$ is finite dimensional and $B_0$ is $p$-nilpotent. 
We can replace $L$ with $L/B$. So we can assume that $A'=0$. 
In particular, $A_1\su M$. We claim that $L$ is solvable. Indeed, let $H=(L_0+M)/A$. It follows from Theorem \ref{Engel's Theorem} that  $H'$ is nilpotent.
Thus, $H$ is solvable and, since $A$ is abelian, we deduce that  $L_0+M$ is solvable. But $L/(L_0+M)$ is abelian which implies  that $L$ is solvable, as claimed. Now, we argue by induction on
the derived length $s$ of $L$.
Suppose first that $L$ is metabelian, that is $(L', L')=0$.  Since $(L_0, L_0)$ is abelian,  it follows that $u^{p^t}=0$, for every $u\in (L_0, L_0)$. So, $(L_0, L_0)$ is $p$-nilpotent. Now we show that $(M, L_1)$ is $p$-nilpotent. 
Let $y_1, \ldots, y_n$ be linearly independent elements in $M$ so that their images span $M/A_1$. 
If $\dim L_1/M=1$, we take $z\in L_1\setminus M$.
In part (1), we proved that $(M, L_1)$ is $p$-nil.
Thus, there exists an integer $m>t$ such that $(y_i, y_j)^{p^m}=(z, y_k)^{p^m}=0$,
for all $1\leq i, j, k\leq n$. Then 
\begin{align*}
(\sum_{i=1}^n \alpha_i y_i, \beta z{+}\sum_{j=1}^n \beta_j y_j)^{p^m}&{=}
\bigg( \sum_{i=1}^n \sum_{j=1}^n \alpha_i\beta_j (y_i, y_j) +\sum_{i=1}^n  \alpha_i\beta (y_i, z)\bigg)^{p^m}\\
&{=}\sum_{i=1}^n \sum_{j=1}^n  (\alpha_i\beta_j y_i, y_j)^{p^m} {+}\sum_{i=1}^n   (\alpha_i\beta y_i, z)^{p^m}=0,
\end{align*}
where $\beta$ and the $\alpha_i$ and $\beta_j$ are in $\F$.
On the other hand, let $y\in L_1$ and $x\in A_1$. Since $x^2=0$, we get 
\begin{align*}
[x, y]^2=(x, y)^2-2(x,y, y)x.
\end{align*}
Since $A'=0$, we deduce that $[x, y]^{2p}=(x, y)^{2p}$.
Thus, $(x, y)^{2p^t}=[x, y]^{2p^t}=0$. 
We deduce that $(M, L_1)^{p^m}=0$. Since $L'$ is abelian, it is clear that $((L_0, L_0)+(M, L_1))^{p^m}=0$.
Now suppose $\delta_s(L)=0$, for some $s\geq 3$.
Let $H=\la \delta_{s-2} (L)\ra_p$ and $K=\la \delta_{s-1} (L)\ra_p$. 
Note that $H_1\su [L_1, L_0]\su M$.
Since $H$ is a metabelian ideal of $L$, we have $((H_0, H_0)+ (H_1, H_1))^{p^m}=0$. 
On the other hand, by the induction hypothesis applied to $L/K$,  we have  $((L_0, L_0)+ (M, L_1))^{p^m}\su  K_0$. But $K_0=\la (H_0, H_0)+ (H_1, H_1)\ra_p$. So, we get $((L_0, L_0)+ (M, L_1))^{p^{2m}}=0$, as required.
\qed

Note that the difference between the restricted case and ordinary case mentioned in the introduction arises from Equation \eqref{imprtant-equation}.
In the ordinary Lie superalgebra case, discussed in \cite{BRU}, Equation \eqref{imprtant-equation} immediately implies that 
$(y,y)=\alpha (y, z)$ and $(z, z)=\alpha^{-1} (y, z)$, for some $\alpha\in \F$. We then deduce that $y$ and $z$ must be dependent and there is no need for the field to be perfect.

\begin{lem}\label{ideal-nil}
Let $L$ be a restricted Lie superalgebra and $A$  a homogeneous ideal of $L$ of finite codimension.
Let $I$ be an ideal of $u(A)$ that is stable under the adjoint action of $L$.
If $u(L)$ is PI and  $I$ is nil of bounded index then so is  $Iu(L)$.
\end{lem}
\Proof Let $R=u(L)$.
Let $x_1,\ldots, x_n\in L_0$ and $y_1, \ldots, y_m\in L_1$ such that
$$
L=A+\la x_1,\ldots, x_n, y_1, \ldots, y_m\ra_{\F}.
$$
We assume that the $x_i$ and $y_j$ are linearly independent modulo $A$.
By the PBW Theorem, $R$ has a basis consisting of the monomials of the form
\begin{align*}
&x_1^{\alpha_1}\ldots x_n^{\alpha_n}y_1^{\beta_1}\ldots y_m^{\beta_m}w\\
&0\leq \alpha_i<p, \quad \beta_j\in \{0,1\},
\end{align*}
where the $w$'s are PBW monomials in $u(A)$.
Let $D=u(A)$. Note that $R$ is a right $D$-module of finite rank $r=p^n2^m$.
Now consider the regular representation 
$$
\rho: R\to \mbox{End}(R_D)=M_r(D),
$$ where 
$\rho(u):R\to R$ is defined by $\rho(u)v=uv$, for every $u, v\in R$. 
Note that $\rho$ is injective because  $R$ is unital.
Thus, under $\rho$, we can embed $R$ into $M_r(D)$.
Since $I$ is stable under the adjoint action of $L$, we have
$RI=IR=RIR$. We claim that $RI$ embeds into $M_r(I)$.
Indeed, let $v_1\in RI$ and $v_2\in R$.
Since $\rho(v_1)(v_2)=v_1v_2\in RI$, we can write $v_1v_2$ as a linear combinations of elements of the form $ua$, where $u\in R$ and $a\in I$. So, each $ua$ and hence $v_1v_2$ is a linear combination of elements of the form $x_1^{\alpha_1}\ldots x_n^{\alpha_n}y_1^{\beta_1}\ldots y_m^{\beta_m}b$, where $b\in I$, as claimed.

Therefore, it suffices to show that $M_r(I)$ is nil of bounded index.
Recall  Levitzki's Theorem and Shirshov's Height Theorem
stating that every $t$-generated PI algebra which is nil of bounded index $s$ is (associative) nilpotent of a bound given as a function of $s$, $t$, and $d$, where $d$ is the degree of the polynomial identity, see \cite{Lev} and \cite{Shirshov}.
So, if $S$ is any $r^2$-generated subalgebra of $I$ then there exists a constant $k$ such that
 $S^k=0$. Now, let $T\in M_r(I)$ and denote by $S$ the subalgebra of $I$ generated by all entries of $T$.
So,  $T^i\in M_r(S^i)$, for every $i$. Since $S^k=0$, we get $T^k=0$. Since $k$ is independent of $T$,
it follows that $M_r(I)$ is nil of bounded index, as required.
\qed

\noindent\emph{Proof of the Main Theorem.}
The implication   $(3)\Rightarrow (1)$ is obvious  while 
$(1)\Rightarrow (2)$  follows from  Proposition \ref{M}. To prove
$(2)\Rightarrow (3)$, we shall use the fact that the class of   nil algebras of bounded index is closed under extensions.
If $\dim L_1/M=1$ then we may assume that $L_1=M+\F z$, where $(z, z)$ is not $p$-nilpotent.
By Theorem \ref{Petrogradsky}, there exists a homogeneous restricted ideal $A$ of finite codimension in $L$ such that
 $A'$ is finite-dimensional and $(A')_0$ is  $p$-nilpotent.
Note that,  by Lemma \ref{L-res-nilpotent}, $(A, A)$ is nilpotent.
Thus, $\w(\la A'\ra_p)$ is associative nilpotent, by Lemma \ref{w(L)-nilpotent}. Hence, $A'u(L)$ is associative  nilpotent.
So, we can replace $L$ with  $L/(A')_p$ to assume
that $A$ is abelian. We claim that $(A_1, L_1)$ is $p$-nilpotent. If $M=L_1$, then 
$(L_1, L_1)$ is $p$-nilpotent and the claim is obvious. So  we may assume that $(z, z)$ is not $p$-nilpotent
and prove  $A_1\su M$. Suppose that there exists $x\in A_1$ such that $x\notin M$.
Since $\dim L/M=1$,  we may assume that $x=z+\alpha y$, for some $y\in M$ and $\alpha\in \F$.
Then, since $(x,x)=0$, we get
$$
(z, z)=(x-\alpha y, x-\alpha y)=-2\alpha(x, y)+\alpha^2 (y, y)=(-2\alpha x+\alpha^2 y, y).
$$
 Since $y\in M$,  it follows from  the hypothesis  that $(y, L_1)$ is $p$-nilpotent. So, $(z, z)$ must be $p$-nilpotent, which is a contradiction.
 Note that  $u(A)$ is the tensor product of a commutative algebra with the Grassmann algebra. Thus, $u(A)$ satisfies 
$[x, y, z]=0$. So, we have $(x+y)^p=x^p+y^p$, for all $x, y\in u(A)$. 
Let $B=\la (A, L)\ra_p$. Note that $B$ is an  abelian ideal of $L$  and  since $(A_0, L_0)$ and $(A_1, L_1)$ are both $p$-nilpotent, $\w(B)$ is nil of bounded index.
It follows that $I=Bu(A)$ is nil of bounded index.
Furthermore, by Lemma \ref{ideal-nil}, $Iu(L)=Bu(L)$ is also a nil ideal of $u(L)$ of bounded index.
But $Iu(L)$ is the kernel of the homomorphism $R\to u(L/B)$. Thus, we can replace
$L$ with $L/B$ to assume that $A$ is central in $L$.
It follows that $L'$ is finite-dimensional. Note that $A_1u(A)$ is  nil of index $p$. Since $A_1$ is a central ideal of $L$, $A_1u(L)$ is nil of bounded index, by Lemma \ref{ideal-nil}. Thus, we may assume that $A_1=0$.
It follows that $M$ is finite-dimensional. Thus, by Lemma \ref{w(L)-nilpotent}, $((M, L_1)+M)u(L)$ is associative  nilpotent. We can now assume $M=0$. Hence, $L_1=\F z$ and $(z, L_0)=0$. Let $H=L_0$. Since $H'$ is finite-dimensional  and $p$-nilpotent, we deduce, by Lemma \ref{w(L)-nilpotent}, that $\w(\la H'\ra_p)$ is associative nilpotent. So we can replace $L$ with $L/\la H'\ra_p$.
Hence, $(L_0, L)=0$. It then follows that  $[R, R]=0$, as required.
\qed

\section{Examples}\label{examples}
We provide examples showing that the restriction on the field in the main result is necessary.

\begin{eg}\label{cntexm}
Let $L=L_0\oplus L_1$ be a restricted Lie superalgebra over a non-perfect field $\F$, where $L_0=\la x_1, x_2, x_3\ra_{\F}$ and $L_1=\la y, z\ra_{\F}$. We assume  that $(L_0, L)=0$ and set $x_1=(y,y)$, $x_2=(z, z)$, and $x_3=(y,z)$. Let $\alpha\in \F$ be an element whose $p$th root does not lie in $\F$.
We define the $p$-mapping by setting $x_1^p=x_1$, $x_2^p=\alpha^2 x_1$, and $x_3^p=\alpha x_1$. The following statements hold:
\begin{enumerate}
\item The commutator ideal of $u(L)$ is nil of bounded index; hence $u(L)$  satisfies a non-matrix PI.
\item $(c, c)$ is not $p$-nilpotent, for every $c\in L_1$.
\end{enumerate}
\end{eg}
\Proof Note that  $[y, z]^2=x_3^2-x_1x_2$ is a central element in $u(L)$ and $[y, z]^{2p}=0$.  Every element in the commutator ideal of $u(L)$ 
is of the form $u=(\alpha y+\beta z+\gamma)[y, z]$, where $\alpha, \beta, \gamma$ are in the center of  $u(L)$.
We observe that 
\begin{align*}
u^2&=(\alpha y+\beta z+\gamma)[y, z](\alpha y+\beta z+\gamma)[y, z]\\
&=(\alpha y+\beta z+\gamma)^2[y, z]^2+ (\alpha y+\beta z+\gamma)[y, z, \alpha y+\beta z+\gamma][y, z]\\
&=(\alpha y+\beta z+\gamma)^2[y, z]^2+ (\alpha y+\beta z+\gamma)( -2\alpha y-2\beta z)[y, z]^2\\
&=(\alpha y+\beta z+\gamma)( -\alpha y-\beta z+\gamma)[y, z]^2
\end{align*} 
Thus $u^{2p}=0$, for every $u\in [u(L), u(L)]u(L)$. Hence, the commutator ideal of $u(L)$ is nil of index $2p$.

Next we prove (2). Suppose to the contrary that there exists $c\in L_1$ such that $(c, c)$ is $p$-nilpotent.
Without loss of generality, we may assume $c=y+\beta z$. We have
$$
[c, z]=-(c,z)+2cz, \quad [c, z]^2=(c, z)^2-(c,c)(z, z).
$$
Since $[c, z]^{2p}=0$, we get   $(c, z)^{2p}=(c, c)^{p}(z, z)^{p}$. Since $(c, c)$ is $p$-nilpotent, 
   $(c, z)$ must be  $p$-nilpotent.
So there exists an integer $m$ such that  $(y+\beta z, z)^{p^m}=0$.
Note that 
\begin{align*}
(y+\beta z, z)^{p^m}=x_3^{p^m}+\beta^{p^m} x_2^{p^m}&=(\alpha x_1)^{p^{m-1}}+\beta^{p^m}(\alpha^2 x_1)^{p^{m-1}}\\
&=\alpha^{p^{m-1}}x_1+\beta^{p^m}\alpha^{2p^{m-1}}x_1=0
\end{align*}
We get $\beta^{p^m}\alpha^{p^{m-1}}+1=0$. Hence, $(\beta^p\alpha+1)^{p^{m-1}}=0$
which implies that $\alpha=(-\beta)^p$. So $\alpha$ must have $p$th root in $\F$ which is a contradiction.
\qed

One might ask if similar examples exist when $\dim L_1\geq 3$. The answer is yes as the following example shows.

\begin{eg}\label{cntexmII} Let $\F\su K$ be a  field extension and suppose that there exist  $u, v\in K$ so that 
 $\{ 1, u, v\}$ is linearly independent over $\F$ and $u^p, v^p\in \F$. For example, one can take 
 $K=\F_p(u, v)$ and $\F=\F_p(u^p, v^p)$, where $\F_p$ is the prime field. Indeed, here we have $[K: \F]=p^2$.

  We define a restricted Lie superalgebra $L=L_0\oplus L_1$ over $\F$ as follows.  
Let $L_1=\la y_1, y_2, y_3\ra_{\F}, L_0=\la z_{ij}, x_k\mid 1\leq i<j\leq 3, 1\leq k\leq 3\ra_{\F}$ and assume that $(L_0, L)=0$. Set $z_{ij}=(y_i, y_j), x_k=(y_k, y_k)$, for all $1\leq i<j\leq 3$ and $1\leq k\leq 3$. 
Let $\alpha=u^p$, $\beta=v^p$ and define the $p$-mapping by 
$x_1^p=x_1$, $x_2^p=\alpha^2 x_1$, $x_3^p=\beta^2 x_1$, 
$z_{12}^p=\alpha x_1$,  $z_{13}^p=\beta x_1$, and $z_{23}^p=\alpha\beta x_1$ . The following statements hold:
\begin{enumerate}
\item The commutator ideal of $u(L)$ is  nil of bounded index.
\item $(c, c)$ is not $p$-nilpotent, for every $c\in L_1$.
\end{enumerate}
\end{eg}
\Proof  Let $J$ be the Jacobson radical of $u(L)$. So, $J$ is associative nilpotent. 
Let $A=L_0\oplus \la y_1, y_2\ra$ and set $P=u(A)$. Then,  by Example \ref{cntexm},   $[P, P]P$ is a nil ideal of $P$ of bounded index. So, by Lemma \ref{ideal-nil}, we deduce that
$[P, P]u(L)$ is nil of bounded index. Hence, $[P, P]u(L)\su J$. We can  similarly prove  that
$[y_1, y_3], [y_2, y_3]\in J$. This implies that the commutator ideal of $u(L)$ is nil of bounded index.
 This proves (1).

To prove (2), suppose to the contrary that there exists $c\in L_1$ such that $(c, c)$ is $p$-nilpotent.
Without loss of generality, we may assume $c=y_3+\alpha_1 y_1 +\alpha_2 y_2$. Then
since the commutator ideal of $u(L)$ is nil of bounded index,   $(y_3+\alpha_1 y_1 +\alpha_2 y_2, y_3)$ is  $p$-nilpotent.
So there exists $m$ such that  $(y_3+\alpha_1 y_1 +\alpha_2 y_2, y_3)^{p^m}=0$.
Note that 

\begin{align*}
(y_3{+}\alpha_1 y_1 {+}\alpha_2 y_2, y_3)^{p^m}&{=}x_3^{p^m}+\alpha_1^{p^m} z_{13}^{p^m}+\alpha_2^{p^m} z_{23}^{p^m}\\
&=(\beta^2 x_1)^{p^{m-1}}+\alpha_1^{p^m} (\beta x_1)^{p^{m-1}}+\alpha_2^{p^m} (\alpha\beta x_1)^{p^{m-1}}\\
&=\beta^{p^{m-1}} x\bigg(   \beta ^{p^{m-1}}+\alpha_1^{p^m} + \alpha_2^{p^m} \alpha^{p^{m-1}}\bigg)=0.
\end{align*}
We get $(\beta+\alpha_1^{p} + \alpha_2^{p} \alpha) ^{p^{m-1}}=0$. Thus, $\beta+\alpha_1^{p} + \alpha_2^{p} \alpha=0$ and so $v+\alpha_1 + \alpha_2u=0$.
But this means that $u, v, 1$ are linearly dependent over $\F$, a contradiction.


\begin{thebibliography}{9}

\bibitem[AK]{AK} E.V. Aladova,  A.N. Krasilnikov,
Polynomial identities in nil-algebras,
\emph{Trans. Amer. Math. Soc.} \textbf{361} (2009), no. 11, 5629--5646.

\bibitem[Am]{Am}
S.A. Amitsur, A generalization of Hilbert's Nullstellensatz,
\emph{Proc. Amer. Math. Soc.} \textbf{8} (1957), 649--656.

\bibitem[BMPZ]{BMPZ} Y. Bahturin, A. Mikhalev, V.M. Petrogradsky, M. Zaicev,
\emph{Infinite Dimensional Lie Superalgebras}, (Walter de Gruyter, Berlin, 1992).

\bibitem[B74]{B74} Y. Bahturin, Identities in the universal envelopes of Lie algebras,
\emph{ J. Austral. Math. Soc.} \textbf{18} (1974), 10--21.

\bibitem[BRU]{BRU} J. Bergen, D.M. Riley, H. Usefi, Lie superalgebras whose  enveloping algebras satisfying a non-matrix polynomial identity, \emph{Israel J. Math}, to appear.

\bibitem[BRT]{BRT} Y. Billig, D. Riley, V. Tasic, Nonmatrix varieties and nil-generated algebras whose units satisfy a group identity \emph{J. Algebra}  \textbf{190}  (1997),  no. 1, 241--252.




\bibitem[K96]{K96} A.R. Kemer, PI-algebras and nil algebras of bounded index,
\emph{Trends in Ring Theory}, CMS Conf. Proc., Miskolc, (1996), (Amer. Math. Soc., Providence, 1998), \textbf{22}, 59--69.

 \bibitem[K91]{Kemer2} A.R. Kemer, \emph{Ideal of Identities of Associative Algebras},
  (Amer. Math. Soc., Providence, RI, 1991), Vol. 87.


\bibitem[L80]{L80} V.N. Latyshev, Nonmatrix varieties of associative algebras,
\emph{Mat. Zametki}  \textbf{27}  (1980), no. 1, 147--156.

\bibitem[L77]{L77} V.N. Latyshev, The complexity of nonmatrix varieties of associative algebras, I, II.
\emph{Algebra i Logika}  \textbf{16}  (1977), no. 2, 149--183, 184--199, 249--250.

\bibitem[L63]{L63}  V.N. Latyshev, Two remarks on $PI$-algebras,
\emph{Sibirsk. Mat. Zh.} \textbf{4} (1963), 1120--1121.

\bibitem[Lev]{Lev} J. Levitzki, On a problem of A. Kurosch,
\emph{Bull. Amer. Math. Soc.} \textbf{52} (1946), 1033--1035.

\bibitem[MPR]{MPR} S.P. Mishchenko, V.M. Petrogradsky,  A. Regev, Characterization of non-matrix varieties of
associative algebras, \emph{Israel J. Math.} to appear.




 \bibitem[P]{P90} D.S. Passman, Enveloping algebras satisfying a polynomial identity,
  \emph{J. Algebra.} \textbf{134(2)} (1990), 469--490.

\bibitem[P92]{P92} V.M. Petrogradski, Identities in the enveloping algebras for modular Lie superalgebras,
\emph{J. Algebra} \textbf{145} (1992), no. 1, 1--21.

\bibitem[P91]{P91}  V.M. Petrogradsky, The existence of an identity in a restricted envelope,
 \emph{Mat. Zametki} \textbf{49(1)} (1991), 84--93.


\bibitem[RW]{RW99} D.M. Riley, M.C. Wilson, Group algebras and enveloping algebras with nonmatrix and semigroup identities,
\emph{Comm. Algebra} \textbf{27} (7) (1999), 3545--3556.

\bibitem[R97]{R97} D.M. Riley, PI-algebras generated by nilpotent elements of bounded index,
\emph{J. Algebra}  \textbf{192 } (1997),  no. 1, 1--13.



\bibitem[Sam]{Sam} L.M. Samoilov, On the radical of a relatively free associative algebra over fields of positive characteristic, \emph{Mat. Sb.} \textbf{199} (2008), no. 5, 81--126.

\bibitem[Sch]{sch} M. Scheunert, The theory of Lie superalgebras,
\emph{Lecture Notes in Math.} \textbf{716} (1979).


\bibitem[Sh]{Shirshov} A.I. Shirshov, On rings with identity relations,
\emph{Mat. Sb.} 43 (85) (1957), 277--283.


\end{thebibliography}
\end{document}